\documentclass[brochure,english,12pt]{bourbaki}
\usepackage{amssymb,amsfonts,amsmath,footnote,comment}
\usepackage[francais]{babel}
\usepackage[all]{xy}
\SelectTips{cm}{}

\date{Juin 2013}
\bbkannee{65\`eme ann\'ee, 2012-2013}
\bbknumero{1072}
\title{Categorification of Lie algebras}
\subtitle{d'apr\`es Rouquier, Khovanov-Lauda, ...}
\author{Joel KAMNITZER}
\address{University of Toronto\\
Department of Mathematics\\
40 St. George Street\\
Toronto, Ontario \\
Canada M5S 2E4}
\email{jkamnitz@math.toronto.edu}

\def\sl{\mathfrak{sl}}
\def\C{\mathbb{C}}

\def\Hom{\operatorname{Hom}}
\def\sA{\mathcal{A}}
\def\dU{\dot{U}}
\def\Z{\mathbb{Z}}
\def\Vect{\mathcal{V}ect}
\def\Cat{\mathcal{C}at}
\def\cU{\mathcal{U}}
\def\g{\mathfrak{g}}
\def\Z{\mathbb{Z}}
\def\N{\mathbb{N}}
\def\Seq{\operatorname{Seq}}
\def\bi{\mathbf{i}}
\def\cC{\mathcal{C}}
\def\cV{\mathcal{V}}
\newcommand{\scs}{\scriptstyle}

\def\pmod{\mathrm{pmod}}
\newcommand{\negspace}{\! \!}

\newcommand{\Uupdot}[1]{
   \xy {\ar (0,-3)*{};(0,3)*{} };(0,0)*{\bullet};(1.5,0)*{};(-1.5,0)*{};(2,-2)*{\scs #1};\endxy}

\newcommand{\Udowndot}[1]{
   \xy {\ar (0,3)*{};(0,-3)*{} };(0,0)*{\bullet};(1.5,0)*{};(-1.5,0)*{};(2,-2)*{\scs #1};\endxy}

\newcommand{\Ucupr}[1]{\;\;
    \vcenter{\xy (-2,2)*{}; (2,2)*{} **\crv{(-2,-2) & (2,-2)}?(1)*\dir{>};
            (2,-3)*{};(4,2)*{\scs #1}; \endxy} \;\; }

\newcommand{\Ucupl}[1]{\;\;
    \vcenter{\xy (2,2)*{}; (-2,2)*{} **\crv{(2,-2) & (-2,-2)}?(1)*\dir{>};
            (2,-3)*{};(3,1.5)*{\scs #1}; \endxy} \;\; }

\newcommand{\Ucapr}[1]{\;\;
    \vcenter{\xy (-2,-1)*{}; (2,-1)*{} **\crv{(-2,3) & (2,3)}?(1)*\dir{>};
            (2,-3)*{};(3,2)*{\scs #1}; \endxy} \;\; }

\newcommand{\Ucapl}[1]{\;\;
    \vcenter{\xy (2,-2)*{}; (-2,-2)*{} **\crv{(2,2) & (-2,2)}?(1)*\dir{>};
            (2,-3)*{};(3,1.5)*{\scs #1}; \endxy} \;\; }

\newcommand{\Ucross}[2]{\;\;
    \vcenter{\xy {\ar (2.5,-2.5)*{};(-2.5,2.5)*{}}; {\ar (-2.5,-2.5)*{};(2.5,2.5)*{}};
    (3,-1)*{\scs #1};(-3,-1)*{\scs #2};\endxy} \;\; }

\newcommand{\Ucrossd}[2]{
    \vcenter{\xy {\ar (2.5,2.5)*{};(-2.5,-2.5)*{}}; {\ar (-2.5,2.5)*{};(2.5,-2.5)*{} }; (3,-1)*{\scs #1};(-3,-1)*{\scs #2};\endxy}  }

\newcommand{\cyclo}{
  \xy
{\ar (-12,-6);(-12,6)};
{\ar (-4,-6);(-4,6)};
 (4,0)*{\cdots};
{\ar (12,-6);(12,6)};
(-12, 0)*{\bullet}; (-15, 0)*{n_{i_1}};
 (-12,-9)*{i_1}; (-4,-9)*{ i_2};(12,-9)*{ i_m };
  \endxy }

\newcommand{\nilHecke}{
\vcenter{\xy 0;/r.17pc/:
    (-4,-4)*{};(4,4)*{} **\crv{(-4,-1) & (4,1)}?(1)*\dir{>};
    (4,-4)*{};(-4,4)*{} **\crv{(4,-1) & (-4,1)}?(1)*\dir{>};
    (-4,4)*{};(4,12)*{} **\crv{(-4,7) & (4,9)}?(1)*\dir{>};
    (4,4)*{};(-4,12)*{} **\crv{(4,7) & (-4,9)}?(1)*\dir{>};
 \endxy} \
  = 0,\quad \quad
 \vcenter{
 \xy 0;/r.17pc/:
    (-4,-4)*{};(4,4)*{} **\crv{(-4,-1) & (4,1)}?(1)*\dir{>};
    (4,-4)*{};(-4,4)*{} **\crv{(4,-1) & (-4,1)}?(1)*\dir{>};
    (4,4)*{};(12,12)*{} **\crv{(4,7) & (12,9)}?(1)*\dir{>};
    (12,4)*{};(4,12)*{} **\crv{(12,7) & (4,9)}?(1)*\dir{>};
    (-4,12)*{};(4,20)*{} **\crv{(-4,15) & (4,17)}?(1)*\dir{>};
    (4,12)*{};(-4,20)*{} **\crv{(4,15) & (-4,17)}?(1)*\dir{>};
    (-4,4)*{}; (-4,12) **\dir{-};
    (12,-4)*{}; (12,4) **\dir{-};
    (12,12)*{}; (12,20) **\dir{-};
\endxy} \;
  = \;
 \vcenter{
 \xy 0;/r.17pc/:
    (4,-4)*{};(-4,4)*{} **\crv{(4,-1) & (-4,1)}?(1)*\dir{>};
    (-4,-4)*{};(4,4)*{} **\crv{(-4,-1) & (4,1)}?(1)*\dir{>};
    (-4,4)*{};(-12,12)*{} **\crv{(-4,7) & (-12,9)}?(1)*\dir{>};
    (-12,4)*{};(-4,12)*{} **\crv{(-12,7) & (-4,9)}?(1)*\dir{>};
    (4,12)*{};(-4,20)*{} **\crv{(4,15) & (-4,17)}?(1)*\dir{>};
    (-4,12)*{};(4,20)*{} **\crv{(-4,15) & (4,17)}?(1)*\dir{>};
    (4,4)*{}; (4,12) **\dir{-};
    (-12,-4)*{}; (-12,4) **\dir{-};
    (-12,12)*{}; (-12,20) **\dir{-};
\endxy} \; , \quad \quad
 \xy 0;/r.18pc/:
  (4,4);(4,-4) **\dir{-}?(0)*\dir{<}+(2.3,0)*{};
  (-4,4);(-4,-4) **\dir{-}?(0)*\dir{<}+(2.3,0)*{};
 \endxy
  = \negspace \negspace
\xy 0;/r.18pc/:
  (0,0)*{\xybox{
    (-4,-4)*{};(4,4)*{} **\crv{(-4,-1) & (4,1)}?(1)*\dir{>}?(.25)*{\bullet};
    (4,-4)*{};(-4,4)*{} **\crv{(4,-1) & (-4,1)}?(1)*\dir{>};
     (-10,0)*{};(10,0)*{};
     }};
  \endxy \negspace \negspace
  - \negspace \negspace
\xy 0;/r.18pc/:
  (0,0)*{\xybox{
    (-4,-4)*{};(4,4)*{} **\crv{(-4,-1) & (4,1)}?(1)*\dir{>}?(.75)*{\bullet};
    (4,-4)*{};(-4,4)*{} **\crv{(4,-1) & (-4,1)}?(1)*\dir{>};
     (-10,0)*{};(10,0)*{};
     }};
  \endxy \negspace \negspace
 = \negspace \negspace
\xy 0;/r.18pc/:
  (0,0)*{\xybox{
    (-4,-4)*{};(4,4)*{} **\crv{(-4,-1) & (4,1)}?(1)*\dir{>};
    (4,-4)*{};(-4,4)*{} **\crv{(4,-1) & (-4,1)}?(1)*\dir{>}?(.75)*{\bullet};
     (-10,0)*{};(10,0)*{};
     }};
  \endxy \negspace \negspace
  - \negspace \negspace
  \xy 0;/r.18pc/:
  (0,0)*{\xybox{
    (-4,-4)*{};(4,4)*{} **\crv{(-4,-1) & (4,1)}?(1)*\dir{>} ;
    (4,-4)*{};(-4,4)*{} **\crv{(4,-1) & (-4,1)}?(1)*\dir{>}?(.25)*{\bullet};
     (-10,0)*{};(10,0)*{};
     }};
  \endxy
}

\newcommand{\KLR}{
 \vcenter{\xy 0;/r.17pc/:
    (-4,-4)*{};(4,4)*{} **\crv{(-4,-1) & (4,1)}?(1)*\dir{>};
    (4,-4)*{};(-4,4)*{} **\crv{(4,-1) & (-4,1)}?(1)*\dir{>};
    (-4,4)*{};(4,12)*{} **\crv{(-4,7) & (4,9)}?(1)*\dir{>};
    (4,4)*{};(-4,12)*{} **\crv{(4,7) & (-4,9)}?(1)*\dir{>};
  (-5,-3)*{\scs i};
     (5.1,-3)*{\scs j};
 \endxy}
 =
 \begin{cases}
     \xy 0;/r.17pc/:
  (3,8);(3,-8) **\dir{-}?(.5)*\dir{<}+(2.3,0)*{};
  (-3,8);(-3,-8) **\dir{-}?(.5)*\dir{<}+(2.3,0)*{};
  (-5,-6)*{\scs i};     (5.1,-6)*{\scs j};
 \endxy \quad \text{if $\langle \alpha_i, \alpha_j)=0$,}\\ \\
  \vcenter{\xy 0;/r.17pc/:
  (3,8);(3,-8) **\dir{-}?(.5)*\dir{<}+(2.3,0)*{};
  (-3,8);(-3,-8) **\dir{-}?(.5)*\dir{<}+(2.3,0)*{};
  (-3,4)*{\bullet};(-6.5,5)*{};
  (-5,-6)*{\scs i};     (5.1,-6)*{\scs j};
 \endxy} \;\; - \;\;
  \vcenter{\xy 0;/r.17pc/:
  (3,8);(3,-8) **\dir{-}?(.5)*\dir{<}+(2.3,0)*{};
  (-3,8);(-3,-8) **\dir{-}?(.5)*\dir{<}+(2.3,0)*{};
  (3,4)*{\bullet};(7,5)*{};
  (-5,-6)*{\scs i};     (5.1,-6)*{\scs j};
 \endxy}
   \quad \text{if $i \leftarrow j$,} \\  \\
   \vcenter{\xy 0;/r.17pc/:
  (3,8);(3,-8) **\dir{-}?(.5)*\dir{<}+(2.3,0)*{};
  (-3,8);(-3,-8) **\dir{-}?(.5)*\dir{<}+(2.3,0)*{};
  (3,4)*{\bullet};(-6.5,5)*{};
  (-5,-6)*{\scs i};     (5.1,-6)*{\scs j};
 \endxy} \;\; - \;\;
  \vcenter{\xy 0;/r.17pc/:
  (3,8);(3,-8) **\dir{-}?(.5)*\dir{<}+(2.3,0)*{};
  (-3,8);(-3,-8) **\dir{-}?(.5)*\dir{<}+(2.3,0)*{};
  (-3,4)*{\bullet};(7,5)*{};
  (-5,-6)*{\scs i};     (5.1,-6)*{\scs j};
 \endxy}
   \quad \text{if $i \rightarrow j$.} \\
 \end{cases}
}

\newcommand{\dotslide}{
\xy 0;/r.18pc/:
  (0,0)*{\xybox{
    (-4,-4)*{};(4,4)*{} **\crv{(-4,-1) & (4,1)}?(1)*\dir{>}?(.75)*{\bullet};
    (4,-4)*{};(-4,4)*{} **\crv{(4,-1) & (-4,1)}?(1)*\dir{>};
    (-5,-3)*{\scs i};
     (5.1,-3)*{\scs j};
     (-10,0)*{};(10,0)*{};
     }};
  \endxy
 \;\; =
\xy 0;/r.18pc/:
  (0,0)*{\xybox{
    (-4,-4)*{};(4,4)*{} **\crv{(-4,-1) & (4,1)}?(1)*\dir{>}?(.25)*{\bullet};
    (4,-4)*{};(-4,4)*{} **\crv{(4,-1) & (-4,1)}?(1)*\dir{>};
    (-5,-3)*{\scs i};
     (5.1,-3)*{\scs j};
     (-10,0)*{};(10,0)*{};
     }};
  \endxy
\qquad  \xy 0;/r.18pc/:
  (0,0)*{\xybox{
    (-4,-4)*{};(4,4)*{} **\crv{(-4,-1) & (4,1)}?(1)*\dir{>};
    (4,-4)*{};(-4,4)*{} **\crv{(4,-1) & (-4,1)}?(1)*\dir{>}?(.75)*{\bullet};
    (-5,-3)*{\scs i};
     (5.1,-3)*{\scs j};
     (-10,0)*{};(10,0)*{};
     }};
  \endxy
\;\;  =
  \xy 0;/r.18pc/:
  (0,0)*{\xybox{
    (-4,-4)*{};(4,4)*{} **\crv{(-4,-1) & (4,1)}?(1)*\dir{>} ;
    (4,-4)*{};(-4,4)*{} **\crv{(4,-1) & (-4,1)}?(1)*\dir{>}?(.25)*{\bullet};
    (-5,-3)*{\scs i};
     (5.1,-3)*{\scs j};
     (-10,0)*{};(12,0)*{};
     }};
  \endxy
}

\newcommand{\Rthreeeasy}{
 \vcenter{
 \xy 0;/r.17pc/:
    (-4,-4)*{};(4,4)*{} **\crv{(-4,-1) & (4,1)}?(1)*\dir{>};
    (4,-4)*{};(-4,4)*{} **\crv{(4,-1) & (-4,1)}?(1)*\dir{>};
    (4,4)*{};(12,12)*{} **\crv{(4,7) & (12,9)}?(1)*\dir{>};
    (12,4)*{};(4,12)*{} **\crv{(12,7) & (4,9)}?(1)*\dir{>};
    (-4,12)*{};(4,20)*{} **\crv{(-4,15) & (4,17)}?(1)*\dir{>};
    (4,12)*{};(-4,20)*{} **\crv{(4,15) & (-4,17)}?(1)*\dir{>};
    (-4,4)*{}; (-4,12) **\dir{-};
    (12,-4)*{}; (12,4) **\dir{-};
    (12,12)*{}; (12,20) **\dir{-};
  (-6,-3)*{\scs i};
  (6,-3)*{\scs j};
  (15,-3)*{\scs k};
\endxy}
 \;\; =\;\;
 \vcenter{
 \xy 0;/r.17pc/:
    (4,-4)*{};(-4,4)*{} **\crv{(4,-1) & (-4,1)}?(1)*\dir{>};
    (-4,-4)*{};(4,4)*{} **\crv{(-4,-1) & (4,1)}?(1)*\dir{>};
    (-4,4)*{};(-12,12)*{} **\crv{(-4,7) & (-12,9)}?(1)*\dir{>};
    (-12,4)*{};(-4,12)*{} **\crv{(-12,7) & (-4,9)}?(1)*\dir{>};
    (4,12)*{};(-4,20)*{} **\crv{(4,15) & (-4,17)}?(1)*\dir{>};
    (-4,12)*{};(4,20)*{} **\crv{(-4,15) & (4,17)}?(1)*\dir{>};
    (4,4)*{}; (4,12) **\dir{-};
    (-12,-4)*{}; (-12,4) **\dir{-};
    (-12,12)*{}; (-12,20) **\dir{-};
  (7,-3)*{\scs k};
  (-6,-3)*{\scs j};
  (-14,-3)*{\scs i};
\endxy}
}

\newcommand{\Rthreehard}{
\;\; \vcenter{
 \xy 0;/r.17pc/:
    (-4,-4)*{};(4,4)*{} **\crv{(-4,-1) & (4,1)}?(1)*\dir{>};
    (4,-4)*{};(-4,4)*{} **\crv{(4,-1) & (-4,1)}?(1)*\dir{>};
    (4,4)*{};(12,12)*{} **\crv{(4,7) & (12,9)}?(1)*\dir{>};
    (12,4)*{};(4,12)*{} **\crv{(12,7) & (4,9)}?(1)*\dir{>};
    (-4,12)*{};(4,20)*{} **\crv{(-4,15) & (4,17)}?(1)*\dir{>};
    (4,12)*{};(-4,20)*{} **\crv{(4,15) & (-4,17)}?(1)*\dir{>};
    (-4,4)*{}; (-4,12) **\dir{-};
    (12,-4)*{}; (12,4) **\dir{-};
    (12,12)*{}; (12,20) **\dir{-};
  (-6,-3)*{\scs i};
  (6,-3)*{\scs j};
  (14,-3)*{\scs i};
\endxy}
\ - \
 \vcenter{
 \xy 0;/r.17pc/:
    (4,-4)*{};(-4,4)*{} **\crv{(4,-1) & (-4,1)}?(1)*\dir{>};
    (-4,-4)*{};(4,4)*{} **\crv{(-4,-1) & (4,1)}?(1)*\dir{>};
    (-4,4)*{};(-12,12)*{} **\crv{(-4,7) & (-12,9)}?(1)*\dir{>};
    (-12,4)*{};(-4,12)*{} **\crv{(-12,7) & (-4,9)}?(1)*\dir{>};
    (4,12)*{};(-4,20)*{} **\crv{(4,15) & (-4,17)}?(1)*\dir{>};
    (-4,12)*{};(4,20)*{} **\crv{(-4,15) & (4,17)}?(1)*\dir{>};
    (4,4)*{}; (4,12) **\dir{-};
    (-12,-4)*{}; (-12,4) **\dir{-};
    (-12,12)*{}; (-12,20) **\dir{-};
  (6,-3)*{\scs i};
  (-6,-3)*{\scs j};
  (-14,-3)*{\scs i};
\endxy}
 \; =\;
\begin{cases}
\xy 0;/r.17pc/:
  (4,12);(4,-12) **\dir{-}?(.5)*\dir{<};
  (-4,12);(-4,-12) **\dir{-}?(.5)*\dir{<}?(.25)*\dir{};
  (12,12);(12,-12) **\dir{-}?(.5)*\dir{<}?(.25)*\dir{};
  (-6,-9)*{\scs i};     (6.1,-9)*{\scs j};
  (14,-9)*{\scs i};
 \endxy \quad \text{if $i \leftarrow j $,} \\ \\
 -\xy 0;/r.17pc/:
  (4,12);(4,-12) **\dir{-}?(.5)*\dir{<};
  (-4,12);(-4,-12) **\dir{-}?(.5)*\dir{<}?(.25)*\dir{};
  (12,12);(12,-12) **\dir{-}?(.5)*\dir{<}?(.25)*\dir{};
  (-6,-9)*{\scs i};     (6.1,-9)*{\scs j};
  (14,-9)*{\scs i};
 \endxy \quad \text{if $i \rightarrow j $.}
 \end{cases}
}

\newcommand{\biadjoint}{
 \xy   0;/r.17pc/:
    (-8,0)*{}="1";
    (0,0)*{}="2";
    (8,0)*{}="3";
    (-8,-10);"1" **\dir{-};
    "1";"2" **\crv{(-8,8) & (0,8)} ?(0)*\dir{>} ?(1)*\dir{>};
    "2";"3" **\crv{(0,-8) & (8,-8)}?(1)*\dir{>};
    "3"; (8,10) **\dir{-};
     \endxy \quad
     = \negspace
\xy   0;/r.17pc/:
    (-8,0)*{}="1";
    (0,0)*{}="2";
    (8,0)*{}="3";
    (0,-10);(0,10)**\dir{-} ?(.5)*\dir{>};
     \endxy \negspace
= \quad
    \xy   0;/r.17pc/:
    (8,0)*{}="1";
    (0,0)*{}="2";
    (-8,0)*{}="3";
    (8,-10);"1" **\dir{-};
    "1";"2" **\crv{(8,8) & (0,8)} ?(0)*\dir{>} ?(1)*\dir{>};
    "2";"3" **\crv{(0,-8) & (-8,-8)}?(1)*\dir{>};
    "3"; (-8,10) **\dir{-};
    \endxy
 \quad \quad \quad \xy  0;/r.17pc/:
    (8,0)*{}="1";
    (0,0)*{}="2";
    (-8,0)*{}="3";
    (8,-10);"1" **\dir{-};
    "1";"2" **\crv{(8,8) & (0,8)} ?(0)*\dir{<} ?(1)*\dir{<};
    "2";"3" **\crv{(0,-8) & (-8,-8)}?(1)*\dir{<};
    "3"; (-8,10) **\dir{-};
    \endxy \quad
     = \negspace
\xy  0;/r.17pc/:
    (8,0)*{}="1";
    (0,0)*{}="2";
    (-8,0)*{}="3";
    (0,-10);(0,10)**\dir{-} ?(.5)*\dir{<};
    \endxy \negspace
     = \quad
\xy   0;/r.17pc/:
    (-8,0)*{}="1";
    (0,0)*{}="2";
    (8,0)*{}="3";
    (-8,-10);"1" **\dir{-};
    "1";"2" **\crv{(-8,8) & (0,8)} ?(0)*\dir{<} ?(1)*\dir{<};
    "2";"3" **\crv{(0,-8) & (8,-8)}?(1)*\dir{<};
    "3"; (8,10) **\dir{-};
    \endxy
}

\newcommand{\invert}{
\left( \xy 0;/r.15pc/:
    (-4,-4)*{};(4,4)*{} **\crv{(-4,-1) & (4,1)}?(1)*\dir{>} ;
    (4,-4)*{};(-4,4)*{} **\crv{(4,-1) & (-4,1)}?(0)*\dir{<};
  \endxy \;\;, \
  \vcenter{\xy 0;/r.15pc/:
 (-4,0)*{}="t1"; (4,0)*{}="t2";
 "t2";"t1" **\crv{(4,7) & (-4,7)}; ?(0)*\dir{<}
 ?(.8)*\dir{};;
 \endxy}\ , \
 \vcenter{\xy 0;/r.15pc/:
 (-4,0)*{}="t1"; (4,0)*{}="t2";
 "t2";"t1" **\crv{(4,7) & (-4,7)}; ?(0)*\dir{<}
 ?(.8)*\dir{}+(0,-.1)*{\bullet}+(-3,2.5)*{\scs 1};;
 \endxy}
  \ , \dots,
 \vcenter{\xy 0;/r.15pc/:
 (-4,0)*{}="t1"; (4,0)*{}="t2";
 "t2";"t1" **\crv{(4,7) & (-4,7)}; ?(0)*\dir{<}
 ?(.8)*\dir{}+(0,-.1)*{\bullet}+(-3,2.5)*{\scs {r-1}};;
 \endxy} \
 \right)
 : E_iF_i \rightarrow F_iE_i \oplus I_\lambda^{\oplus r}
 }

\newcommand{\EiFj}{
 \vcenter{   \xy 0;/r.18pc/:
    (-4,-4)*{};(4,4)*{} **\crv{(-4,-1) & (4,1)}?(1)*\dir{>};
    (4,-4)*{};(-4,4)*{} **\crv{(4,-1) & (-4,1)}?(1)*\dir{<};?(0)*\dir{<};
    (-4,4)*{};(4,12)*{} **\crv{(-4,7) & (4,9)};
    (4,4)*{};(-4,12)*{} **\crv{(4,7) & (-4,9)}?(1)*\dir{>};
 (-6,-3)*{\scs i};
     (6,-3)*{\scs j};
 \endxy}
 \;\;= \;\;
\xy 0;/r.18pc/:
  (3,9);(3,-9) **\dir{-}?(.55)*\dir{>}+(2.3,0)*{};
  (-3,9);(-3,-9) **\dir{-}?(.5)*\dir{<}+(2.3,0)*{};
  (-5,-6)*{\scs i};     (5.1,-6)*{\scs j};
 \endxy
 \qquad
    \vcenter{\xy 0;/r.18pc/:
    (-4,-4)*{};(4,4)*{} **\crv{(-4,-1) & (4,1)}?(1)*\dir{<};?(0)*\dir{<};
    (4,-4)*{};(-4,4)*{} **\crv{(4,-1) & (-4,1)}?(1)*\dir{>};
    (-4,4)*{};(4,12)*{} **\crv{(-4,7) & (4,9)}?(1)*\dir{>};
    (4,4)*{};(-4,12)*{} **\crv{(4,7) & (-4,9)};
  (-6,-3)*{\scs i};
     (6,-3)*{\scs j};
 \endxy}
 \;\;=\;\;
\xy 0;/r.18pc/:
  (3,9);(3,-9) **\dir{-}?(.5)*\dir{<}+(2.3,0)*{};
  (-3,9);(-3,-9) **\dir{-}?(.55)*\dir{>}+(2.3,0)*{};
  (-5,-6)*{\scs i};     (5.1,-6)*{\scs j};
 \endxy
}

\newcommand{\newcross}{
\xy
  (0,0)*{\xybox{
    (-4,-4)*{};(4,4)*{} **\crv{(-4,-1) & (4,1)}?(1)*\dir{>} ;
    (4,-4)*{};(-4,4)*{} **\crv{(4,-1) & (-4,1)}?(0)*\dir{<};
    (-5,-3)*{\scs i};
     (-5,3)*{\scs j};
     (-12,0)*{};(12,0)*{};
     }};
  \endxy \negspace \negspace
:=
 \xy 0;/r.19pc/:
  (0,0)*{\xybox{
    (4,-4)*{};(-4,4)*{} **\crv{(4,-1) & (-4,1)}?(1)*\dir{>};
    (-4,-4)*{};(4,4)*{} **\crv{(-4,-1) & (4,1)};
     (-4,4);(-4,12) **\dir{-};
     (-12,-4);(-12,12) **\dir{-};
     (4,-4);(4,-12) **\dir{-};(12,4);(12,-12) **\dir{-};
     (-10,0)*{};(10,0)*{};
     (-4,-4)*{};(-12,-4)*{} **\crv{(-4,-10) & (-12,-10)}?(1)*\dir{<}?(0)*\dir{<};
      (4,4)*{};(12,4)*{} **\crv{(4,10) & (12,10)}?(1)*\dir{>}?(0)*\dir{>};
      (-14,11)*{\scs j};(-2,11)*{\scs i};
      (14,-11)*{\scs j};(2,-11)*{\scs i};
     }};
  \endxy
  \ = \
  \xy 0;/r.19pc/:
  (0,0)*{\xybox{
    (-4,-4)*{};(4,4)*{} **\crv{(-4,-1) & (4,1)}?(1)*\dir{<};
    (4,-4)*{};(-4,4)*{} **\crv{(4,-1) & (-4,1)};
     (4,4);(4,12) **\dir{-};
     (12,-4);(12,12) **\dir{-};
     (-4,-4);(-4,-12) **\dir{-};(-12,4);(-12,-12) **\dir{-};
     (10,0)*{};(-10,0)*{};
     (4,-4)*{};(12,-4)*{} **\crv{(4,-10) & (12,-10)}?(1)*\dir{>}?(0)*\dir{>};
      (-4,4)*{};(-12,4)*{} **\crv{(-4,10) & (-12,10)}?(1)*\dir{<}?(0)*\dir{<};
      (14,11)*{\scs i};(2,11)*{\scs j};
      (-14,-11)*{\scs i};(-2,-11)*{\scs j};
     }};
  \endxy
}

\begin{document}
\maketitle

\noindent{\bf INTRODUCTION}

Categorification is the process of finding hidden higher level structure.  To categorify a natural number, we look for a vector space whose dimension is that number.  For example, the passage from Betti numbers to homology groups was an important advance in algebraic topology.

To categorify a vector space $ V $, we look for a category $\cC $ whose Grothendieck group is that vector space, $K(\cC) = V $.  If $V $ carries an action of a Lie algebra $ \g $, then it is natural to look for functors $ F_a : \cC \rightarrow \cC $ for each generator $ a $ of $ \g $, such that $F_a $ gives the action of $ a $ on the Grothendieck group level.  In this case, we say that we have categorified the representation $ V $.

There are two general motivations for trying to categorify representations.  First, by studying the category $ \cC $, we hope to learn more about the vector space $ V$.  For example, we get a special basis for $ V $ coming from classes of indecomposable objects of $ \cC $.  Second, we may use the action of $ \g $ on $\cC $ to learn more about $ \cC $.  For example, Chuang-Rouquier used categorification to prove Brou\'e's abelian defect group conjecture for symmetric groups.

Recently, there has been amazing progress towards constructing categorifications of representations of semisimple (or more generally Kac-Moody) Lie algebras.  In this report, we aim to give an introduction to this theory.  We start with the categorification of $ \sl_2 $ and its representations.  We explain the naive definition and then the ``true'' definition, due to Chuang-Rouquier \cite{CR}.  We also explain how this definition leads to interesting equivalences of categories.    We then address general Kac-Moody Lie algebras, reaching the definition of the Khovanov-Lauda-Rouquier 2-category \cite{R1,KL3}.  We explain the relationship to Lusztig's categories of perverse sheaves, due to Varagnolo-Vasserot \cite{VV} and Rouquier \cite{R2}.  We close by discussing three fundamental examples of categorical representations: modular representation theory of symmetric groups (due to Lascoux-Leclerc-Thibon \cite{LLT}, Grojnowski \cite{Gr}, and Chuang-Rouquier \cite{CR}), cyclotomic quotients of KLR algebras (due to Kang-Kashiwara \cite{KK} and Webster \cite{W}), and quantized quiver varieties (due to Zheng \cite{Z} and Rouquier \cite{R2}).

In order to keep the exposition readable, we have made a number of simplifications and glossed over many details.  In particular, we only address simply-laced Kac-Moody Lie algebras (and when it comes to the geometry, only finite-type).  We suggest that interested readers consult the literature for more details.

Throughout this paper, we work over $ \C $; all vector spaces are $ \C$-vector spaces (sometimes they are actually $\C(q)$-vector spaces) and all additive categories are $ \C$-linear.

I would like to thank R. Rouquier, M. Khovanov, and A. Lauda for developing the beautiful mathematics which is presented here and for their many patient explanations (an extra thank you to A. Lauda for allowing me to use his diagrams).  I also thank D. Ben-Zvi, R. Bezrukavnikov, A. Braverman, J. Brundan, C. Dodd, D. Gaitsgory, H. Nakajima, A. Kleshchev, A. Licata, D. Nadler, B. Webster, G. Williamson, and O. Yacobi for interesting discussions about categorification over many years and a special thank you to S. Cautis for our long and fruitful collaboration.  Finally, I thank S. Cautis, M. Khovanov, A. Lauda, C. Liu, S. Morgan, R. Rouquier, B. Webster and O. Yacobi for their helpful comments on a first draft of this paper.

\bigskip

\section{Categorification of $\mathfrak{sl}_2 $ representations}
\subsection{The structure of finite-dimensional representations} \label{se:repsl2}
The Lie algebra $ \sl_2(\C) $ has the basis $$ e = \left[ \begin{smallmatrix} 0 & 1 \\ 0 & 0 \end{smallmatrix} \right], \ h = \left[ \begin{smallmatrix} 1 & 0 \\ 0 & -1 \end{smallmatrix} \right], \ f = \left[ \begin{smallmatrix}0 & 0 \\ 1 & 0 \end{smallmatrix} \right].$$
Consider a finite-dimensional representation $ V $ of $ \sl_2 $.  A basic theorem of representation theory states that $ h $ acts semisimply on $ V $ with integer eigenvalues.  Thus we may write $ V = \oplus_{r \in \mathbb Z} V_r $ as the direct sum of the eigenspaces for $ h $.  Moreover the commutation relations between the generators $ e, f, h $ imply the following.
\begin{enumerate}
\item For each $ r$, $ e $ restricts to a linear map $ e  : V_r \rightarrow V_{r+2}$.
\item Similarly, $ f $ restricts to a linear map $ f : V_r \rightarrow V_{r-2} $.
\item These restrictions obey the commutation relation
\begin{equation} \label{eq:effe}
 ef - fe|_{V_r} = r I_{V_r}.
 \end{equation}
\end{enumerate}
Conversely, a graded vector space $ V = \oplus V_r $, along with raising and lowering operators $ e, f $ as above, defines a representation of $ \sl_2 $ if these operators satisfy the relation \eqref{eq:effe}.

The following example will be very instructive.
\begin{exem} \label{eg:finitesets}
Let $ X $ be a finite set of size $ n $.  Let $ V = \C^{P(X)} $ be a vector space whose basis consists of the subsets of $ X $.  For $ r = -n, -n + 2, \dots, n $, define $ V_r $ to be the span of subsets of size $ k $, where $ r = 2k-n $.

Define linear maps $ e : V_r \rightarrow V_{r+2}$, $ f : V_r \rightarrow V_{r-2} $ by the formulas
\begin{equation} \label{eq:EFsubset}
e(S) = \sum_{T \supset S, |T| = |S| + 1} T, \quad \quad f(S) = \sum_{T \subset S, |T| = |S| - 1} T
\end{equation}
It is easy to check that $ (ef - fe)(S) = (2k-n) S $, if $ S $ has size $k$. (The basic reason is that there are $ n-k $ ways to add something to $ S $ and $ k $ ways to take something away from $ S $.)

Thus this defines a representation of $ \sl_2$.  In fact, this representation is isomorphic to an $n$-fold tensor product $ (\C^2)^{\otimes n} $ of the standard representation of $ \sl_2$.
\end{exem}

We will also need the concept of a representation of the quantum group $ U_q \sl_2 $, though we will neither need nor give an explicit definition of $ U_q \sl_2$.

For each integer $ r $, let
$$ [r] := \frac{q^r - q^{-r}}{q - q^{-1}} = q^{r-1} + q^{r-3} + \cdots + q^{-r+1}
$$
denote the quantum integer (the second expression is only valid if $ r \ge 0 $).

A representation of $ U_q \sl_2 $ is a graded $ \C(q) $ vector space $ V = \oplus V_r $  along with raising $ e : V_r \rightarrow V_{r+2} $ and lowering $ f : V_r \rightarrow V_{r-2} $ operators such that $ ef - fe|_{V_r} = [r] I_{V_r} $.

\subsection{Naive categorical action}
Once we think of an $\sl_2 $ representation in terms of a sequence of vector spaces together with raising and lowering operators, we are led to the notion of an action of $ \sl_2 $ on a category.
\begin{defi}
A naive categorical $\sl_2 $ action consists of a sequence $ D_r $ of additive categories along with additive functors $ E : D_r \rightarrow D_{r+2} $, $ F : D_r \rightarrow D_{r-2} $, for each $ r $, such that there exist isomorphisms of functors
\begin{align}
EF|_{D_r} &\cong FE|_{D_r} \oplus I_{D_r}^{\oplus r}, \ \text{ if } r \ge 0 \label{eq:EFFE1}\\
FE|_{D_r} &\cong EF|_{D_r} \oplus I_{D_r}^{\oplus r}, \ \text{ if } r \le 0 \label{eq:EFFE2}
\end{align}
\end{defi}

Suppose that the categories $ D_r $ carry a naive categorical $ \sl_2$ action.  Then we can construct a usual $ \sl_2 $ representation as follows.  We set $ V_r = K(D_r) $, the complexified split Grothendieck group.  The functors $ E, F $ give rise to linear maps $ e : V_r \rightarrow V_{r+2} $, $ f : V_r \rightarrow V_{r-2} $ and we can easily see that \eqref{eq:EFFE1} and \eqref{eq:EFFE2} give the commutation relation \eqref{eq:effe}.  Thus we get a representation of $ \sl_2 $ on $ V = \oplus V_r $.  We say that the categories $ D_r $ categorify the representation $ V = \oplus V_r $.

It is also useful to consider a graded version of the above definition.  A graded additive category is a category $ \cC $ along with an additive functor $ \langle 1 \rangle: \cC \rightarrow \cC $.  We define a graded naive categorical $ \sl_2 $ action as above but with \eqref{eq:EFFE1}, \eqref{eq:EFFE2} replaced by
\begin{align*}
EF|_{D_r} &\cong FE|_{D_r} \oplus I_{D_r}\langle r-1 \rangle \oplus \cdots \oplus I_{D_r}\langle -r+1 \rangle, \ \text{ if } r \ge 0 \\
FE|_{D_r} &\cong EF|_{D_r} \oplus I_{D_r}\langle r-1 \rangle \oplus \cdots \oplus I_{D_r}\langle -r+1 \rangle, \ \text{ if } r \le 0
\end{align*}
The Grothendieck groups $ K(D_r) $ will then carry an action of $ U_q \sl_ 2$.

We will now give an example of a naive categorical action which will build on Example \ref{eg:finitesets}.

In Example \ref{eg:finitesets}, we studied subsets of a finite set.  There is a well-known analogy between subsets of an $n$-element set and subspaces of an $n$-dimensional vector space over a finite field $\mathbb F_q$, where $ q $ is a power of a prime.  This analogy suggests that we try to construct a representation of $ \sl_2 $ on $ \oplus V_r$, where $ V_r = \C^{G(k, \mathbb F_q^n)} $ is a $ \C $-vector space whose basis is $ G(k,\mathbb F_q^n) $, the set of $k$-dimensional subspaces of $ \mathbb F_q^n $ (where $ r = 2k-n $ as before).  If we define $ e,f $ as in (\ref{eq:EFsubset}), then we get a representation of the quantum group $U_{\sqrt{q}} \sl_2$ (after a slight modification).

The finite set $ G(k,\mathbb F_q^n) $ is the set of $ \mathbb F_q $-points of a projective variety, called the Grassmannian.  By Grothendieck's fonctions-faisceaux correspondence, we can categorify $ \C^{G(k, \mathbb F_q^n)} $ using an appropriate category of sheaves on $ G(k, \mathbb{\overline{F}}_q^n) $.  For simplicity, we switch to characteristic 0 and consider sheaves on $ G(k,\C^n) $, the Grassmannian of $k$-dimensional subspaces of $ \C^n $.

For each $ r = -n, -n+2, \dots, n$, we let $ D_r = D_c^b(G(k,\C^n)) $ denote the bounded derived category of constructible sheaves (again here $ r = 2k-n$).  These are graded categories, where the grading comes from homological shift.  With the above motivations, we will define a categorical $ \sl_2 $ action using these categories.

For each $ k $, we define the 3-step partial flag variety
$$ Fl(k,k+1, \C^n) = \{ 0 \subset V \subset V' \subset \C^n : \dim V = k, \dim V' = k+1 \} $$
$Fl(k,k+1, \C^n) $ serves as a correspondence between $ G(k,\C^n) $ and $ G(k+1, \C^n) $ and thus it can be used to define functors between categories of sheaves on these varieties.  Let $ p : Fl(k,k+1, \C^n) \rightarrow G(k,\C^n)$ and $ q : Fl(k,k+1, \C^n) \rightarrow G(k+1,\C^n)$ denote the two projections.

We define
\begin{align*}
E : D_r = D_c^b(G(k,\C^n)) &\rightarrow D_{r+2} = D_c^b(G(k+1,\C^n)) \\
\sA &\mapsto q_*(p^* \sA)
\end{align*}
\begin{align*}
F : D_r &\rightarrow D_{r-2} \\
\sA &\mapsto p_*(q^* \sA)
\end{align*}
The above definition of $ E, F $ parallels the definition (\ref{eq:EFsubset}).

The following result was proven in an algebraic context (i.e. after applying the Bei\-lin\-son-Bernstein correspondence) by Bernstein-Frenkel-Khovanov \cite{BFK}.
\begin{theo} \label{th:naivesl2}
This defines a graded naive categorical $ \sl_2$ action.
\end{theo}
The proof of this theorem is relatively straightforward.  To illustrate the idea, let us fix $ V \in G(k,\C^n) $ and consider $ A_1 = \{ V' : V \subset V', \dim V' = k+1 \} $ and $ A_2 = \{V' : V \supset V', \dim V' = k-1 \}$; these are the varieties of ways to increase or decrease $ V $.  Note that $ A_1 $ is a projective space of dimension $ n-k-1$ and $ A_2 $ is a projective space of dimension $ k-1$.  Thus $ \dim H^*(A_2) - \dim H^*(A_1) = 2k - n$.  This observation combined with the decomposition theorem proves the above result.

\begin{rema}
The Grothendieck group of these categories $ D_r $ is actually infinite-dimensional.  To cut down to a finite dimensional situation, we can consider the full subcategories $ D'_r = P_{Sch}(G(k,\C^n)) $ consisting of direct sums of homological shifts of IC-sheaves on Schubert varieties.  The subcategories $ D'_r $ carry a naive categorical $ \sl_2 $ action and by considering dimensions of weight spaces, we can see that they categorify the representation $ (\C^2)^{\otimes n} $.
\end{rema}

\subsection{Categorical $\sl_2$-action}
In the definition of naive categorical $ \sl_2 $ action, we only demanded that there exist isomorphisms of functors in (\ref{eq:EFFE1}) and (\ref{eq:EFFE2}). We did not specify the data of these isomorphisms.  This is very unnatural from the point of view of category theory.  However, it is not immediately obvious how to specify these isomorphisms nor what relations these isomorphisms should satisfy.

In their breakthrough paper, Chuang-Rouquier \cite{CR} solved this problem.  First, it is natural to assume that the functors $ E, F $ be adjoint (this is a categorification of the fact that $e, f $ are adjoint with respect to the Shapovalov form on any finite-dimensional representation of $ \sl_2 $).

Now (assume $ r \ge 0 $), we desire to specify a isomorphism of functors
$$
(\phi, \psi_0, \dots, \psi_{r-1}) : EF|_{D_r} \rightarrow FE|_{D_r} \oplus I_{D_r}^{\oplus r}
$$
so $ \phi \in \Hom(EF,FE) \cong \Hom(EE,EE)$ (using the adjunction) and $ \psi_s \in \Hom(EF, I) \cong \Hom(E,E) $ (again using the adjunction).  Thus it is natural to choose two elements $ T \in \Hom(EE,EE) $ and $ X \in \Hom(E,E) $ such that $ \phi $ corresponds to $ T $ and $ \psi_s $ corresponds to $ X^s $ for $ s = 0 ,\dots, r-1$.

This leads us to the following definition, essentially due to Chuang-Rouquier \cite{CR}.
\begin{defi} \label{def:catsl2action}
A categorical $ \sl_2$ action consists of
\begin{enumerate}
\item a sequence $ D_r $ of additive categories, with $ D_r = 0 $ for $ r \ll 0 $,
\item functors $ E : D_r \rightarrow D_{r+2} $, $ F : D_r \rightarrow D_{r-2} $, for each $ r $,
\item natural transformations $ \varepsilon : EF \rightarrow I, \ \eta : I \rightarrow FE, \ X : E \rightarrow E, \ T : E^2 \rightarrow E^2 $
\end{enumerate}
such that the following holds.
\begin{enumerate}
\item The morphisms $ \varepsilon, \eta $ are the units and counits of adjunctions.
\item If $ r \ge 0 $, the morphism
\begin{align}
(\sigma, \varepsilon, \varepsilon \circ XI_F\dots, \varepsilon\circ X^{r-1}I_F) : EF|_{D_r} \rightarrow FE|_{D_r} \oplus I_{D_r}^{\oplus r}
\end{align}
is an isomorphism, where $ \sigma : EF \rightarrow FE $ is defined as the composition $$ EF \xrightarrow{\eta I_{EF}} FEEF \xrightarrow{I_F T 1_F} FEEF \xrightarrow{I_{FE} \varepsilon} FE. $$ (And we impose a similar isomorphism condition if $ r \le 0 $.)
\item The morphisms $ X, T $ obey the following relations.
\begin{enumerate}
\item In $ \Hom(E^2, E^2)$, we have $ XI_E \circ T - T \circ I_E X = I_{E^2} = T \circ XI_E - I_E X \circ T $.
\item In $ \Hom(E^2, E^2)$, we have $T^2 = 0$.
\item In $\Hom(E^3, E^3)$, we have $ T I_E \circ I_E T \circ T I_E = I_E T \circ T I_E \circ 1_E T$.
\end{enumerate}
\end{enumerate}
\end{defi}

\begin{rema}
If we work in the graded setting, then it is natural to ask that $ X $ have degree 2, i.e. that it be a morphism $ X: E \rightarrow E\langle 2\rangle $.  Likewise, we give $ T $ degree $ -2 $.  The degrees of $ \varepsilon $ and $ \eta $ depend on $ r $.
\end{rema}

At first glance, it is not apparent where the relations among the $ X, T $ come from.  To motivate them, we introduce the nil affine Hecke algebra.

\begin{defi}
The nil affine Hecke algebra $ H_n $ is the algebra with generators $ x_1, \dots, x_n, t_1, \dots, t_{n-1} $ and relations
\begin{gather*}
t_i^2 = 0,  \ t_i t_{i+1} t_i = t_{i+1} t_i t_{i+1}, \ t_i t_j = t_j t_i  \text{ if }|i-j| > 1, \\
x_i x_j = x_j x_i, \  t_i x_i - x_{i+1} t_i = 1 = x_i t_i - t_i x_{i+1}
\end{gather*}
\end{defi}
Suppose that we have a categorical $ \sl_2 $-action.  Then the morphisms $ X, T $ generate an action of $ H_n $ on $ E^n $.  More precisely, we have an algebra morphism $ H_n \rightarrow \Hom(E^n, E^n) $ by sending $ x_i $ to $ I_{E^{i-1}} X I_{E^{n-i}} $ and $ t_i $ to $ I_{E^{i-1}} T I_{E^{n-i-1}} $.  The above relations among $ X, T $ ensure that the relations of $ H_n $ hold.

\begin{rema}
In their original paper, Chuang-Rouquier \cite{CR} used relations among $ X, T $ modelled after the affine Hecke algebra or degenerate affine Hecke algebra, rather than the nil affine Hecke algebra.  The nil affine Hecke relations were first introduced by Lauda \cite{L}.
\end{rema}

The nil affine Hecke algebra arises quite naturally in the study of the topology of the flag variety.  Let $ Fl(\C^n) $ denote the variety of complete flags in $ \C^n $.  The following result appears to be due to Arabia \cite{A} (see also \cite[Prop. 12.8]{G}).
\begin{prop} \label{prop:nHequalsH}
There is an isomorphism of algebras $$ H_n \cong H_*^{GL_n}(Fl(\C^n) \times Fl(\C^n))$$ where the right hand side carries an algebra structure by convolution.
\end{prop}

\subsection{Categorical $ \sl_2$ actions coming from Grassmannians}
Let us return to constructible sheaves on Grassmannians.  Consider the functor $ E^p : D_c^b(G(k,\C^n)) \rightarrow D_c^b(G(k+p,\C^n))$.  It is given by the correspondence with the partial flag variety
$$
Fl(k,k+1, \dots, k+p, \C^n) = \{ 0 \subset V_0 \subset V_{1} \subset \cdots \subset V_{p} \subset \C^n : \dim V_j = k+ j \}
$$
The map $ Fl(k,k+1, \dots, k+p, \C^n) \rightarrow G(k,\C^n) \times G(k+p,\C^n) $ is a fibre bundle onto its image $ Fl(k, k+p, \C^n) $ with fibre $Fl(\C^p) $.  By Proposition \ref{prop:nHequalsH} this provides an action of the algebra $ H_p $ on the functor $ E^p $.  This can be used to upgrade Theorem \ref{th:naivesl2} to the following result.
\begin{theo} \label{th:sl2actionGrass}
The naive graded categorical $ \sl_2 $ action on $ D_r = D_c^b(G(k,\C^n)) $ extends to a graded categorical $ \sl_2 $ action.
\end{theo}
The above result is well-known but does not appear explicitly in the literature.  It is a special case of the main result of \cite{W2}.

It is worth mentioning a more ``elementary'' version of this categorical $ \sl_2 $ action.  For each $ k = 0, \dots, n$, let $ D''_r $ be the category of finite-dimensional $ H^*(G(k,n)) $-modules (with $ r = 2k-n$).  We have a functor $ D_r \rightarrow D''_r$ given by global sections.  The following result was sketched by Chuang-Rouquier \cite[section 7.7.2]{CR} and a complete proof was given by Lauda \cite[Theorem 7.12]{L}.

\begin{theo}
There exists a categorical $ \sl_2 $ action on $ D''_r $ compatible with the functor $ D_r \rightarrow D''_r $.  This categorifies the $n+1$-dimensional irreducible representation of $ \sl_2 $.
\end{theo}

Moreover, this categorical $\sl_2$ representation is the simplest possible categorification of this irreducible representation; more precisely, it is a minimal categorification, according to the results of Chuang-Rouquier \cite{CR}.

A related construction was given by Cautis, Licata, and the author in \cite{CKL}.  We considered derived categories of coherent sheaves on cotangent bundles to
Grassmannians $ D'''_r := D^b Coh(T^* G(k,\C^n)) $, where again $ r = 2k-n$.  We proved the following result.
\begin{theo} \label{th:coh}
There is a graded categorical $ \sl_2 $ action on $ D'''_r $ where the functors $ E, F $ come from the conormal bundles to the correspondences $ Fl(k, k+1, \C^n) $.  This categorifies the representation $ (\C^2)^{\otimes n} $.
\end{theo}

\subsection{Equivalences}
We will now see how a categorical $ \sl_2 $ action can be used to produce interesting equivalences of categories, following Chuang-Rouquier \cite{CR}.

To motivate the construction, suppose that $ V = \oplus V_r $ is a finite-dimensional representation of $ \sl_2 $.  Then the group $ SL_2 $ acts on $ \oplus V_r $.  In particular the matrix $ s = \left[ \begin{smallmatrix} 0 & 1 \\ -1 & 0  \end{smallmatrix} \right] $ acts on $ V $.  Since $ s $ is a lift of the non-trivial element in the Weyl group of $ SL_2 $, it gives an isomorphism of vector spaces $ s : V_r \rightarrow V_{-r} $ for all $ r $.  We would like to do something similar for categorical $ \sl_2 $ actions.

To do this, let us fix $ r \ge 0 $ and note that the action of $ s $ on $ V_r $ is given by $$ s|_{V_r} = F^{(r)} - E F^{(r+1)} + E^{(2)} F^{(r+2)} - \cdots $$ where $ E^{(n)} = \frac{1}{n!} E^n$.  (Note that this sum is finite since for large enough $ p $, $ V_{r - 2p} = 0 $.)

The alternating signs in this expression suggest that we try to categorify $s $ using a complex.  This complex was introduced by Chuang-Rouquier \cite{CR}, inspired by certain complexes of Rickard.  The following result is due to Chuang-Rouquier \cite{CR} in the abelian case and Cautis-Kamnitzer-Licata \cite{CKL} in the triangulated case (which is the one we state below).
\begin{theo} \label{th:equiv}
Suppose that $ D_r $ is a sequence of triangulated categories carrying a graded categorical $ \sl_2 $ action such that all functors $ E, F $ are exact. Then the complex
$$
S = [F^{(r)} \rightarrow E F^{(r+1)}\langle -1\rangle \rightarrow E^{(2)} F^{(r+2)}\langle -2 \rangle \rightarrow \cdots]
$$
provides an equivalence $ S : D_r \rightarrow D_{-r} $.
\end{theo}
Here $ E^{(n)} $ is defined using a splitting $E^n = {E^{(n)}}^{\oplus n!} $ which is achieved using the action of $ H_n $ on $ E^n $ (see section 4.1.1 of \cite{R1} or section 9.2 of \cite{L}).  The maps in this complex come from the adjunctions.  See section 6.1 of \cite{CR} for more details.

\begin{exem}
Suppose that we have a categorical $ \sl_2 $ action with just $ D_2, D_0, D_{-2} $ non-zero.  Then choosing $ r = 0 $, the above complex has two terms $ S = [I \rightarrow EF \langle -1 \rangle] $. In this case, the equivalence $ S $ is actually a Seidel-Thomas \cite{ST} spherical twist with respect to the functor $ E : D_{-2} \rightarrow D_0 $.  Thus we see that the equivalences coming from categorical $ \sl_2 $ actions generalize the theory of spherical twists.
\end{exem}

Chuang-Rouquier applied Theorem \ref{th:equiv} to prove that certain blocks of modular representations of symmetric groups were derived equivalent.  This proved Brou\'e's abelian defect group conjecture for symmetric groups.  See Theorem \ref{th:actionrepSn} for the construction of the relevant categorical action.

Another very interesting application of Theorem \ref{th:equiv} concerns constructible sheaves on Grassmannians, as in Theorem \ref{th:sl2actionGrass}.  In this case, it can be shown that the resulting equivalence $ D^b_c(G(k,\C^n)) \rightarrow D^b_c(G(n-k, \C^n)) $ is given by the Radon transform.  More precisely, $S $ is given by the integral transform with respect to the kernel $ j_* \C_U $, where $U \subset G(k,\C^n) \times G(n-k, \C^n)$ is the open $ GL_n $-orbit consisting pairs of transverse subspaces\footnote{This result will appear in a forthcoming paper by Cautis, Dodd, and the author}.

Yet another application of Theorem \ref{th:equiv} involves coherent sheaves on cotangent bundles of Grassmannians.  In \cite{CKL}, by combining Theorem \ref{th:equiv} with Theorem \ref{th:coh}, we were able to construct an equivalence $$ D^b Coh(T^* G(k,\C^n)) \rightarrow D^b Coh(T^* G(n-k, \C^n)), $$ thus answering an open problem posed by Kawamata and Namikawa.  (This approach was previously suggested by Rouquier in \cite{RBour}.)  The exact description of the equivalence in this case was given by Cautis \cite{C}.

\section{The Khovanov-Lauda-Rouquier categorification}
We will now rephrase the notion of categorical $ \sl_2 $ action (Definition \ref{def:catsl2action}) from a more general viewpoint.  We will then proceed to define the categorification of any simply-laced Kac-Moody Lie algebra.

\subsection{Generalities on categorification}
Let $ \cC $ be an additive category.  Let  $K(\cC) $ denote the (complexified) split Grothendieck group of $ \cC $; this is the vector space spanned by isomorphism classes $ [A] $ of objects of $ \cC $ modulo the relation $ [A \oplus B] = [A] + [B] $.  If $ \cC $ is a graded additive category, then $ K(\cC) $ is a $ \C[q,q^{-1}]$-module, where we define $ q[A] = [A\langle 1 \rangle] $.  We can then tensor to obtain a $ \C(q)$-vector space, which we will also denote by $ K(\cC) $.

Let $ V $ be a vector space.  A categorification of $ V $ is an additive category $ \cC $, along with an isomorphism of vector spaces $ K(\cC) \cong V $.  If $ V $ is a $ \C(q) $-vector space, then a categorification of $ V $ is a graded additive category $ \cC $, along with an isomorphism of $ \C(q)$-vector spaces $ K(\cC) \cong V $.

We will also need the notion of categorification of algebras.  A monoidal category is an additive category $ \mathcal{C} $, along with an additive bifunctor $ \otimes : \mathcal{C} \times \mathcal{C} \rightarrow \mathcal{C} $, such that $ A \otimes (B \otimes C) = (A \otimes B) \otimes C $\footnote{Actually, this defines the notion of strict monoidal category.}.  If $ \cC $ is a monoidal category, then $ K(\cC) $ acquires the structure of an algebra where the multiplication is defined by $ [A][B] = [A \otimes B] $.

Let $ A $ be an algebra.  A categorification of $ A $ is a monoidal category $ \cC $, along with an isomorphism of algebras $ K(\cC) \cong A $.  (This generalizes in an obvious way to $ \C(q)$-algebras and graded monoidal categories.)

\begin{exem}
The simplest algebra is $ A = \C $.  This algebra is categorified by $ \Vect $, the category of finite-dimensional vector spaces.  Similarly, $ \C(q) $ is categorified by the category of graded vector spaces.

More generally, if $ G $ is a finite group, then the category $Rep(G) $ of finite-dimensional representations of $ G $ categorifies the algebra $ C_c(G) $ of class functions on $ G $.  The isomorphism $ K(Rep(G)) \rightarrow C_c(G) $ is provided by the character map.
\end{exem}

An algebra $ A $ can be regarded as a linear category with one object whose set of endomorphisms is $ A $ and where the composition of morphisms is the multiplication in $ A $.  From this perspective, it is natural to try to categorify more general categories, especially those with very few objects.  To this end, we will need to look at 2-categories.

A 2-category $ \cC $ (for our purposes) is a category enriched over the category of additive categories.  That means we have a set of objects $ \cC $ and for any two objects $ A, B \in \cC $, a category $ \Hom(A, B) $.  We also have associative composition functors $ \Hom(B, C) \times \Hom(A, B) \rightarrow \Hom(A, C) $.  Note that a monoidal category is the same as a 2-category with one object.

The simplest example of a 2-category is $ \Cat $, the 2-category of additive categories.  The objects of $ \Cat $ are additive categories and for any two additive categories $ A, B $, we define $ \Hom(A,B) $ to be the category of functors from $ A $ to $ B $ (the morphisms in $ \Hom(A,B) $ are natural transformations of functors).

If $ \cC $ is a 2-category, then we will define $ K(\cC) $ to be the category whose objects are the same as $ \cC $ and whose morphism sets are defined by $ \Hom_{K(\cC)}(A, B) = K(\Hom(A,B))$.

Let $ A$ be a linear category.  A categorification of $ A $ is an additive 2-category $ \cC $ along with an isomorphism $ K(\cC) \cong A $.

We will also need the notion of idempotent completion (or Karoubi envelope). Recall that if $ \cC $ is an additive category, an idempotent in $ \cC $ is a morphism $ T : A \rightarrow A $ in $ \cC $ such that $ T^2 = T $.  We say that $ T $ splits if we can write $ A $ as a direct sum $ A = A_0 \oplus A_1 $, such that $ T $ acts by $ 0 $ on $ A_0 $ and by $ 1 $ on $ A_1 $.  The idempotent completion $ (\cC)^i $ of $ \cC $ is the smallest enlargement of $ \cC $ such that all idempotents split in $ (\cC)^i $.  If $ \cC $ is a 2-category, then $ (\cC)^i $ will denote the 2-category with the same objects, but where we perform idempotent completion on every Hom-category.

\subsection{2-categorical rephrasing for $ \sl_2 $}

Let us apply this setup to $ A = U\sl_2 $, the universal enveloping algebra.  Actually we will need Lusztig's idempotent form $ \dU \sl_2 $.
Since $ \dU \sl_2 $ carries a system of idempotents, we can regard it as a category.

\begin{defi} \label{def:dUsl2}
The category $ \dU \sl_2 $ has objects $ r \in \Z $.  It is the $ \C$-linear category with generating morphisms $ e \in \Hom(r, r+2) $ and $ f \in \Hom(r, r-2) $, for all $ r $, subject to the relation $ ef - fe = r I_r $ for all $ r $ (this is an equation in $ \Hom(r,r) $).
\end{defi}

A representation of an algebra $ A$ is the same thing as a linear functor $ \dot{A} \rightarrow \Vect $, where $ \dot{A} $ is the category with one object constructed using $ A $.  Thus we can speak more generally of a representation of a linear category $ \cC $ as a linear functor $ \cC \rightarrow \Vect $.

In particular, we can consider linear functors $ \dU \sl_2 \rightarrow \Vect $.   From our discussion in section \ref{se:repsl2},  we can see that a finite-dimensional representation $ V = \oplus V_r $ of $ \sl_2 $ is the same thing as a linear functor $ \dU \sl_2 \rightarrow \Vect $ which takes the object $ r $ to the vector space $ V_r $.

We also have $ \dU_q \sl_2 $, which is defined in the same fashion, except that it is $\C(q) $-linear and the relation is  $ ef - fe = [r] I_r $.

Now we proceed to the question of trying to categorify $ \dU \sl_2 $.  Since it is a category with objects $ r \in \Z $, it will be categorified by a 2-category with the same set of objects.  In the previous section we explained Chuang-Rouquier's definition (Definition \ref{def:catsl2action}) of a categorical $ \sl_2 $ action.  By thinking about this definition, we reach the definition of a 2-category which categorifies $ \dU \sl_2 $.
\begin{defi} \label{def:2catsl2}
Let $\cU \sl_2 $ denote the additive 2-category with
\begin{enumerate}
\item objects $ r\in \Z $,
\item 1-morphisms generated under direct sum and composition by $ E \in \Hom(r,r+2) $ and $ F \in \Hom(r, r-2) $ for all $ r $,
\item 2-morphisms generated by $$ X : E \rightarrow E, \ T : E^2 \rightarrow E^2,\  \eta : I \rightarrow FE,\  \varepsilon: EF \rightarrow I $$
\end{enumerate}
subject to the relations
\begin{enumerate}
\item in $ \Hom(E,E) $, we have $ \varepsilon I_E \circ I_E \eta = I_E $,
\item in $ \Hom(E^2, E^2)$, we have $ XI_E \circ T - T \circ I_E X = I_{E^2} = T \circ XI_E - I_E X \circ T $,
\item in $ \Hom(E^2, E^2)$, we have $T^2 = 0$,
\item in $\Hom(E^3, E^3)$, we have $ T I_E \circ I_E T \circ T I_E = I_E T \circ T I_E \circ 1_E T$,
\item if $ r \ge 0 $, the following 2-morphism
\begin{align}
(\sigma, \varepsilon, \varepsilon \circ X I_F \dots, \varepsilon\circ X^{r-1}I_F) : EF \rightarrow FE \oplus I_r^{\oplus r}
\end{align}
is an isomorphism, where $ \sigma $ is defined as in Definition \ref{def:catsl2action} (plus a similar condition if $ r \le 0 $).
\end{enumerate}
\end{defi}
More precisely, the last condition means that for each $ r$, in the category $\Hom(r,r) $ we adjoin the inverse of $(\sigma, \varepsilon, \varepsilon \circ X I_F \dots, \varepsilon\circ X^{r-1}I_F)$.

Now that we have defined the 2-category $ \cU \sl_2 $, it is natural to consider 2-functors $\cU \sl_2 \rightarrow \Cat $ (these are 2-representations of $ \cU \sl_2 $).  With the above definition, it is easy to see that a categorical $ \sl_2 $ action on some categories $ D_r $ (Definition \ref{def:catsl2action}) is the same thing as a 2-functor $ \cU \sl_2 \rightarrow \Cat $ which takes $ r $ to $ D_r $ for all $ r $.

\begin{rema} \label{rem:compare}
In this definition, we are following Rouquier's definition \cite{R1} of the 2-category.  In the Lauda \cite{L} version, which we denote by $ \cU^L \sl_2 $, we do not invert $(\sigma, \varepsilon, \dots, \varepsilon\circ X^{r-1}I_F)$, but rather add extra relations to ensure that this map is invertible.  In a recent paper, Cautis-Lauda \cite{CL} proved that under some mild assumptions a 2-functor from $\cU \sl_2 $ to $ \Cat $ gives rise to a 2-functor from $\cU^L \sl_2$ to $ \Cat$ (the converse is automatically true).
\end{rema}

The following result is due to Lauda \cite{L}.
\begin{theo} \label{th:sl2cat}
The 2-category $ \cU^L \sl_2 $ categorifies $ \dU \sl_2 $.
\end{theo}

\begin{rema}
The graded version of $\cU \sl_2 $ categorifies Lusztig's $ \dU_q \sl_2 $.  There is also a more precise version of Theorem \ref{th:sl2cat}, which states that the idempotent completion $(\cU^L \sl_2)^i $ categorifies Lusztig's $ \Z[q,q^{-1}] $-form of $ \dU_q \sl_2 $ (if we look at the $\Z[q,q^{-1}]$ version of the Grothendieck group).
\end{rema}

\subsection{The 2-category for general $ \g$}
Suppose that $ \g $ is an arbitrary Kac-Moody Lie algebra.  It is natural to try to extend the above construction from $ \sl_2 $ to $ \g $, in particular to construct a 2-category $ \cU \g $ which categorifies $ \dU \g $.  Roughly equivalent constructions of this 2-category were achieved independently and simultaneously by Khovanov-Lauda \cite{KL1, KL2, KL3} and by Rouquier \cite{R1}.

For simplicity, we will assume that $ \g $ is simply-laced.  Let us fix notation as follows.  Let $ X $ denote the weight lattice of $ \g $.  Let $ I $ denote the indexing set for the simple roots and let $ \alpha_i$ for $ i \in I $ denote the simple roots.  Let $ \Z I \subset X $ be the root lattice and let $ \N I $ denote the positive root cone.  Let $ \langle, \rangle $ denote the symmetric bilinear form on $ X $.  Then $ \langle \alpha_i, \alpha_j \rangle $ are the entries of the Cartan matrix of $ \g $  (these lie in the set $ \{2, -1 ,0 \}$ by assumption).  We choose an orientation of the Dynkin diagram of $ \g $ in order to produce a directed graph, called a quiver and denoted $ Q $.  We write $ i \rightarrow j $ if there is an oriented edge from $ i $ to $ j $ in $ Q $.

The category $ \dU \g $ is constructed from Lusztig's idempotent form of the universal enveloping algebra $ U \g $ and its definition parallels $ \dU \sl_2 $ (Definition \ref{def:dUsl2}).  In particular, it has objects $ \lambda \in X $ and generating morphisms $ e_i \in \Hom(\lambda, \lambda + \alpha_i) $ and $ f_i \in \Hom(\lambda, \lambda - \alpha_i) $ for $ i \in I $ and $ \lambda \in X $ (for reasons of brevity, we do not give a complete list of the relations in $ \dU \g$).  As before, there is a quantum version $ \dU_q \g $ which is obtained by replacing all integers in the definition of $ \dU \g $ by quantum integers.

We will describe the 2-category $ \cU \g $ using graphical notation due to Khovanov and Lauda.  In this graphical notation, 2-morphisms are viewed as string diagrams in the plane, with strings oriented and labelled from $ i \in I $.  The orientations and labels on the strands tell you the source and target of the 2-morphism.  An arrow labelled $ i$ pointing up (resp. down) denotes $ E_i$ (resp. $F_i$).  For more information on this graphical notation see \cite[Section 4]{L}.

\begin{defi} \label{defUcat}
The 2-category $\cU \g $ is defined as follows.
\begin{itemize}
\item The objects are $\lambda$ for $\lambda \in X$.
\item The 1-morphisms are generated by
$$ E_i \in \Hom(\lambda, \lambda + \alpha_i), \quad F_i \in \Hom(\lambda, \lambda-\alpha_i) $$
for $i \in I$ and $\lambda \in X$.
\item The 2-morphisms are generated by
\begin{gather*}
X_i = \Uupdot{i} : E_i \rightarrow E_i, \quad X_i = \Udowndot{i} : F_i \rightarrow F_i, \\ T_{ij} = \Ucross{j}{i} \negspace : E_i E_j \rightarrow E_j E_i, \quad T_{ij} = \Ucrossd{j}{i} : F_i F_j \rightarrow F_j F_i \\
\Ucapr{i}  \negspace : E_i F_i \rightarrow I, \quad \Ucapl{i} \negspace : F_i E_i \rightarrow I, \quad \Ucupr{i} \negspace : I \rightarrow F_iE_i, \quad \Ucupl{i} \negspace : I \rightarrow E_iF_i
 \end{gather*}
 for $ i \in I $ and $ \lambda \in X $.  (We have suppressed $ \lambda $ in the above notation --- it should label a region in each elementary string diagram.  This label tells you the source and target of the $ E_i, F_i $.)
 \end{itemize}
The 2-morphisms are subject to the following relations.
\begin{itemize}
\item The KLR algebra relations among upward pointing string diagrams
\begin{enumerate}
\item
If all strands are labeled by the same $i \in I$, the nil affine Hecke algebra relations
\begin{equation}
 \nilHecke \label{eq_nil_rels}
  \end{equation}

\item For $i \neq j$
\begin{equation}
 \KLR \label{eq_r2_ij-gen}
\end{equation}

\item For $i \neq j$ the dot sliding relations
\begin{eqnarray} \label{eq_dot_slide_ij-gen}
\dotslide
\end{eqnarray}

\item Unless $i = k$ and $ \langle \alpha_i, \alpha_j \rangle < 0$ the relation
\begin{equation}
\Rthreeeasy \label{eq_r3_easy-gen}
\end{equation}
Otherwise, $\langle \alpha_i, \alpha_j \rangle < 0$ and
\begin{equation}
\Rthreehard
 \label{eq_r3_hard-gen}
\end{equation}
\end{enumerate}
\item The cap and cup morphisms are biadjunctions
\begin{equation} \label{eq_biadjoint1}
\biadjoint
\end{equation}

Moreover the dots and crossing are compatible with these biadjunctions.

\item For each $i \ne j $, we have
\begin{equation}
\EiFj \label{eq:EiFj}
\end{equation}
where we define
\begin{equation}
\newcross
\end{equation}
(The equality comes from the biadjointness of the crossing.)

\item For each $ i$ and each $ \lambda $ such that $ r = \langle \lambda, \alpha_i \rangle \ge 0 $, the following 2-morphism is invertible,
\begin{equation}\label{eq:EiFi}
\invert
\end{equation}
Here a dot with a positive integer $k $ indicates that we put $ k$ dots on that strand (in other words, it means $ X_i^k $).  We also impose a similar condition if $ r \le 0 $.
\end{itemize}
\end{defi}

This definition is quite complicated, so let us see where these relations come from.

When $ \g = \sl_2 $, this definition gives the 2-category from Definition \ref{def:2catsl2}.  In fact, \eqref{eq_nil_rels} is relations 2,3,4 from Definition \ref{def:2catsl2} written in diagrammatic form and \eqref{eq_biadjoint1} and \eqref{eq:EiFi} correspond to relations 1 and 5 from Definition \ref{def:2catsl2}. (Actually there is a slight difference, in that the above definition imposes biadjointness, whereas Definition \ref{def:2catsl2} only involves one-sided adjointness.  For more discussion on this see \cite[Theorem 5.16]{R1}.)

Khovanov-Lauda and Rouquier discovered the relations \eqref{eq_r2_ij-gen},\eqref{eq_dot_slide_ij-gen}, \eqref{eq_r3_easy-gen},  and \eqref{eq_r3_hard-gen} based on computations involving cohomology of partial flag varieties and quiver varieties (essentially to get Theorem \ref{th:VV} to hold).

\begin{rema}
As in Remark \ref{rem:compare}, this is the Rouquier version of the 2-category, because of \eqref{eq:EiFi}.  Khovanov-Lauda's version, denoted $\cU^{KL} \g$, bears the same relationship to $ \cU \g$ as Lauda's version, $\cU^L \sl_2$, did in the $ \sl_2 $ case.
\end{rema}

Consider the Grothendieck group $ K(\cU \g) $ as a 1-category.  The generating morphisms are $ e_i = [E_i] , f_i = [F_i] $ as above.  From \eqref{eq:EiFj}, we see that we have $ e_i f_j = f_j e_i $, and from \eqref{eq:EiFi}, we see that $ e_i f_i - f_i e_i = \langle \lambda, \alpha_i \rangle I_\lambda $ in $\Hom(\lambda, \lambda) $.  As these are most of the relations of $ \dU \g $ (there remain the Serre relations), this suggests the following result, which was proven by Khovanov-Lauda \cite{KL3} in the case of $ \sl_n $ and for general $ \g $ by Webster \cite{W}.
\begin{theo}
There is an isomorphism of categories $ K(\cU^{KL} \g) \cong \dU g $.  In other words, the 2-category $ \cU^{KL} \g $ is a categorification of $ \dU \g $.
\end{theo}

\begin{rema} \label{rem:2catgraded}
There is a graded version of $ \cU \g $ with the degree of $ X_i $ equal to $ 2 $ and the degree of $ T_{ij} $ equal to $ - \langle \alpha_i, \alpha_j \rangle $.  This graded version categorifies $ \dU_q \g $.  Again, there is a more precise form relating the idempotent completion of $ \cU^{KL} \g $ and Lusztig's $ \Z[q,q^{-1}]$-form of $ \dU_q \g $.
\end{rema}

\subsection{Categorification of the upper half} \label{se:upperhalf}
It is important to isolate the categorification of the upper half of the envelopping algebra $ U^+ \g $, where $ U^+ \g \subset U \g $ is the subalgebra generated by all $ E_i $ (or equivalently, it is the envelopping algebra of $ \mathfrak{n} $).  Note that $ U^+ \g $ has no idempotents, so we regard it as an algebra, not as a category.  We have the usual grading $ U^+ \g = \oplus_{\nu \in \N I} (U^+ \g)_\nu $.

\begin{defi}
Let $ \cU^+ \g $ denote the monoidal category whose objects are generated (under direct sum and tensor product) by $ E_i $, for $ i \in I $, and whose morphisms are upward pointing string diagrams as in Definition \ref{defUcat} (so the morphisms are generated by the upward pointing dot and crossing with the KLR algebra relations).
\end{defi}

\begin{rema}
In the above definition of $ \cU^+ \g $, the $ E_i $ do not have a source and target object as they do in $ \cU \g $.  Thus $\cU^+ \g $ does not sit inside $ \cU \g $ in any way.  This is not that surprising, as $ U^+ \g $ is not a subalgebra of $ \dU \g $.
\end{rema}

Let $ \nu \in \N I $.  Let $$ \Seq_\nu = \{\bi=(i_1, \dots, i_m) : \alpha_{i_1} + \dots + \alpha_{i_m} = \nu \} $$
This is the set of all ways to write $\nu$ as an ordered sum of simple roots.  For $ \bi \in \Seq_\nu $, we let $ E_\bi = E_{i_1} \cdots E_{i_m} $ (here the $\otimes$-operation in $ \cU^+ \g $ is written as concatenation).  Let  $(\cU^+ \g)_\nu $ denote the full subcategory of $ \cU^+ \g $ whose objects are directs sums of the $ E_\bi $ for $ \bi \in \Seq_\nu $.

We define algebras $ R_\nu := \oplus_{\bi, \bi' \in \Seq_\nu} \Hom_{\cU^+ \g}(E_\bi, E_\bi') $.

These algebras have become known as Khovanov-Lauda-Rouquier (KLR) algebras, though the term quiver Hecke algebras has also been used.  See \cite{B} for a survey paper on these algebras.
\begin{exem}
Suppose that $ \g = \sl_2 $.   Then $ \nu = n \alpha $ for some $n $, $ \Seq_\nu $ has only one element and $ R_\nu = H_n $, the nil affine Hecke algebra.
\end{exem}

By general principles, we have an equivalence of categories $ (\cU^+ \g)_\nu^i \cong R_\nu$-$\pmod $ between the idempotent completion of $ (\cU^+ \g)_\nu $ and the category of projective modules over the KLR algebra $ R_\nu $.  In particular, $ K((\cU^+ \g)_\nu) $ acquires a basis of indecomposable projective $ R_\nu $ modules (these are the same as the indecomposable objects of $ (\cU^+ \g)_\nu $ under the above equivalence).  Note that under the above equivalence, the monoidal structure on $ (\cU^+ \g)^i $ comes from the inclusion $ R_\mu \otimes R_\nu \rightarrow R_{\mu + \nu} $ given by horizontal concatenation of string diagrams.

The following result is due to Khovanov-Lauda \cite[Theorem 1.1]{KL1} (in the simply-laced case).
\begin{theo} \label{th:Upluscat}
$(\cU^+ \g)^i = \oplus_\nu R_\nu$-$\pmod$ is a categorification of $ U^+ \g $.
\end{theo}

In fact, $ K((\cU^+ \g)^i) $ can be given the structure of a bialgebra and then the above result can be strengthened to an isomorphism of bialgebras.

\begin{rema}
There is a graded version of $\cU^+ \g $ which categorifies $ U^+_q \g $.
\end{rema}

\section{Lusztig's perverse sheaves and KLR algebras}
We will now explain a geometric incarnation of the KLR algebras and of the category $ \cU^+ \g$.  For simplicity, let us assume that $ \g $ is of finite type.

\subsection{Lusztig's perverse sheaves}
Recall that we chose an orientation of the Dynkin diagram of $ \g $ to produce the quiver $ Q $.

\begin{defi}
A representation of $ Q $ is a graded vector space $ V = \oplus_{i \in I} V_i $ along with linear maps $ A_{ij} : V_i \rightarrow V_j $ for every directed edge $ i \rightarrow j $ in $ Q $.
\end{defi}
The dimension-vector of a representation $ V $ is defined by $$ \dim V = \sum_{i \in I} \dim V_i \, \alpha_i \in \N I.$$

Let $ M_\nu $ denote the moduli stack of representations of $ Q $ of dimension-vector $ \nu $.  More explicitly, we can present $ M_\nu $ as a global quotient
$$
M_\nu = \bigoplus_{i \rightarrow j} \Hom(\C^{\nu_i}, \C^{\nu_j}) / \prod_i GL(\C^{\nu_i})
$$
where $ \nu = \sum_i \nu_i \alpha_i$.

\begin{exem}
When $ \g = \sl_2$, then the quiver $ Q $ consists of just one vertex with no arrows.  Thus a representation of $ Q $ is just a vector space.  So we see that $ M_{n\alpha} = pt/GL_n$.

When $ \g = \sl_3$, then the quiver $ Q $ consists of two vertices with an arrow between them.  Thus a representation of $ Q $ is a pair of vertices and a linear map between them.  Thus, we see that $ M_{n\alpha_1 + m \alpha_2} = \Hom(\C^n, \C^m)/GL_n \times GL_m $.
\end{exem}

We let $ D(M) := \oplus_\nu D(M_\nu) $ denote the derived category of constructible sheaves on the stack $ M = \sqcup M_\nu $.  Note that we may consider $D(M_\nu) $ as the $ \prod_i GL(\C^{\nu_i})$-equivariant derived category of $\oplus_{(i,j) \in Q} \Hom(\C^{\nu_i}, \C^{\nu_j})$.

Following Lusztig \cite{Lu1,Lu2}, we define a monoidal structure on the category $ D(M) $. We consider the moduli stack of short exact sequences
$$
S = \{ 0 \rightarrow V_1 \rightarrow V_3 \rightarrow V_2 \rightarrow 0 \}
$$
of representations of $ Q $.  We have three projection morphisms $ \pi_1, \pi_2, \pi_3 : S \rightarrow M $ and thus for $ A, B \in D(M) $, we can define $ A * B = {\pi_3}_*(\pi_1^* A \otimes \pi_2^* B) $.  If $ A \in D(M_\nu) $ and $ B\in D(M_\mu) $ then $ A * B \in D(M_{\nu + \mu})$.

The simple perverse sheaves in $ D(M_\nu) $ are precisely the IC-sheaves of the $ \prod_i GL(\C^{\nu_i})$-orbits in
$\oplus_{i \rightarrow j} \Hom(\C^{\nu_i}, \C^{\nu_j})$ and thus are in bijection with the isomorphism classes of representations of $ Q $ of dimension-vector $ \nu $.  Ringel's theorem tells us that the indecomposable representations of $ Q $ have the positive roots as their dimension-vectors.  Thus, the number of isomorphism classes of representations of $ Q $ of dimension-vector $ \nu $ equals the dimension of $ (U^+ \g)_\nu $.

Let $P(M_\nu) $ be the subcategory of $ D(M_\nu) $ consisting of direct sums of homological shifts of simple perverse sheaves in $D(M_\nu)$.  By the decomposition theorem, $ P(M) = \oplus P(M_\nu) $ is a monoidal subcategory.  Note that $ P(M) $ has a graded structure given by homological shift.

Lusztig \cite{Lu1,Lu2} proved the following theorem concerning $ P(M)$.
\begin{theo}
The Grothendieck ring of $ P(M) $ is isomorphic to $ U^+_q \g $.  In other words, $ P(M) $ is a categorification of $ U^+_q \g $.
\end{theo}
By this theorem, $ U^+_q \g$ acquires a basis coming from the classes of the IC-sheaves in $ P(M)$.  This basis is called Lusztig's canonical basis.

\subsection{Relationship to KLR algebras}
It is natural to expect that Lusztig's categorification of $ U^+_q \g $ is related to the categorification of $U^+_q \g $ defined by generators and relations in section \ref{se:upperhalf}.  This result was proven independently by Varagnolo-Vasserot \cite{VV} and by Rouquier \cite{R2}.

\begin{theo} \label{th:VV}
There is an equivalence of additive monoidal categories $ (\cU^+_q)^i \rightarrow P(M) $ defined on generators as follows
\begin{align*}
E_i &\mapsto \C_{M_{\alpha_i}} \\
X_i &\mapsto x_i \in Ext^*(\C_{M_{\alpha_i}}, \C_{M_{\alpha_i}}) \cong H_*^{\C^\times}(pt) \cong \C[x_i] \\
T_{ij} &\mapsto t_{ij}
\end{align*}
\end{theo}
Here we use that $ M_{\alpha_i} \cong pt/\C^\times $.  The definition of $ t_{ij} $ is a bit involved and depends on cases, so we skip the definition.

We can reformulate this theorem using a convolution algebra defined using $ M_\nu $.  We define $ \widetilde{M}_\nu $ to be the moduli stack of complete flags $ 0 \subset V_1 \subset \dots \subset V_m $ of representations of $Q $ with $ \dim V_m = \nu $.  Then we can form the stack $ Z_\nu := \widetilde{M}_\nu \times_{M_\nu} \widetilde{M}_\nu $.  Then $ H_*(Z_\nu) $ is an algebra under convolution.  By Ginzburg \cite[Prop 5.1]{G}, $ H_*(Z_\nu) $ is an Ext-algebra in $ P(M)$.  With this setup, Theorem \ref{th:VV} is equivalent to the existence of compatible isomorphisms $ R_\nu \cong H_*(Z_\nu) $ for all $ \nu $.

\begin{exem}
If we take $ \g = \sl_2 $ and $ \nu = n \alpha $, then $ M_\nu = Fl(\C^n)/GL_n $ and $ Z_\nu = (Fl(\C^n) \times Fl(\C^n))/GL_n $.  In this case the isomorphism $ R_\nu \cong H_*(Z_\nu) $ is precisely the statement of Proposition \ref{prop:nHequalsH}.
\end{exem}

As a corollary of Theorem \ref{th:VV}, we obtain the following result.
\begin{coro}
The basis of $ U^+_q \g $ provided by indecomposable graded projective $ R_\nu$-modules under Theorem \ref{th:Upluscat} is Lusztig's canonical basis.
\end{coro}

\section{Examples of categorical $\g $-representations}

\subsection{Definition of categorical $\g$-representations}
Using the 2-category $ \cU \g $, we can now define a categorical $ \g $-representation to be an additive linear 2-functor $ \cU \g \rightarrow \Cat $.  In particular, a categorical $ \g $-representation consists of a collection of categories $ D_\mu $ for $ \mu \in X $, biadjoint functors $ E_i, F_i : D_\mu \rightarrow D_{\mu \pm \alpha_i} $ and natural transformations $ X_i : E_i \rightarrow E_i,\, T_{ij} : E_i E_j \rightarrow E_j E_i $ satisfying the relations in $ \cU \g$.

A graded categorical $ \g $ representation involves the same setup except that each category $ D_\mu $ is graded, some shifts appear in the biadjointness of $ E_i, F_i $, and the natural transformations $ X_i, T_{ij} $ have degrees as indicated in Remark \ref{rem:2catgraded}.

\subsection{Modular representation theory of symmetric groups}
Going back to the work of Lascoux-Leclerc-Thibon \cite{LLT} and Grojnowski \cite{Gr}, the prime motivating example of a categorical $\g$-representation concerns modular representations of symmetric groups.  In fact, this categorical action has proved to be very important in understanding modular representation theory.

Fix a prime $ p $ and an algebraically closed field $ k $ of characteristic $ p$.  We will be interested in the category $Rep(S_n) $ of finite-dimensional representations of $ S_n $ over $ k $.  These categories will provide an action of the affine Lie algebra $ \widehat{\sl}_p $.  The basic functors we consider between these categories are the induction and restriction functors corresponding to the natural embedding $ S_{n-1} \hookrightarrow S_n $.

Recall the Young-Jucys-Murphy elements $ Y_m := (1 \ m) + \dots + (m-1 \ m) \in kS_n $.  A fundamental result is that the eigenvalues of $ Y_m $ acting on a representation $ M $ lie in the prime subfield $ \Z/p \subset k $.  We will identify of $ \Z/p $ as the set $ I$  of simple roots of our Kac-Moody algebra $\widehat{\sl}_p$.

Let $ i \in \Z/p $.  Using the Young-Jucys-Murphy elements, we define functors $E_i$ and $ F_i $ of $i$-restriction and $i$-induction as follows.  If $ M \in Rep(S_n)$, we let $ E_i(M) $ denote the generalized $i$-eigenspace of $ Y_n $.  Since $ Y_n $ commutes with the action $ S_{n-1} $, we see that $ E_i(M) $ is an $ S_{n-1} $ representation.  Similarly, we define $ F_i(M) $ to be the generalized $i$-eigenspace of $ Y_{n+1} $ acting on $ Ind_{S_n}^{S_{n+1}} M $.

Symmetric polynomials in $ Y_1, \dots, Y_n $ span the centre of $ kS_n $.  Thus we may regard a central character $ \gamma : kS_n \rightarrow k $ as a element of $(\Z/p)^n/S_n $, which we think of as the set of $ n$-element multisubsets of $ \Z/p $.  Thus for each $\mu = (\mu_0, \dots, \mu_{p-1}) \in \N^p $ such that $ \sum \mu_i = n $, we can consider the category $ Rep(S_n)_\mu $ of representations $ M $ of $S_n $ whose generalized central character is given by the multiset $ \gamma(\mu) := \{0^{\mu_0}, \dots, (p-1)^{\mu_{p-1}} \}$.

\begin{theo} \label{th:actionrepSn}
The category $ \oplus_n Rep(S_n) $ carries a categorical $ \widehat{\sl}_p$-action.  More precisely, we get the categorical $ \widehat{\sl}_p$-action as follows
 \begin{itemize}
\item We define $ D_{\omega_0 - \sum \mu_i \alpha_i} = Rep(S_n)_\mu $ for each $ \mu \in \N^p $ (where $ n = \sum \mu_i $).
\item For each $ i \in \Z/p $, we define $ E_i, F_i $ as above.
\item The ``dot'' $ X_i $ and crossing $ T_{ij} $ are defined with the help of $ Y_n $ and the transposition $ (n-1 \ n) $.
\end{itemize}
\end{theo}
The fact that these categories carry a naive categorial $ \widehat{\sl}_p $-action was proven by Lascoux-Leclerc-Thibon \cite{LLT} and by Grojnowski \cite{Gr}. The above statement of an actual categorical $ \widehat{\sl}_p $-action was proven by Chuang-Rouquier \cite[Theorem 4.23]{R2}.  In this theorem, we work with a version of $ \cU \widehat{\sl}_p $ defined over the field $ k $ (rather than $ \C$).

Let us be more precise about the definitions of  $X_i $ and $ T_{ij} $.  Consider the functor $ Res^{S_n}_{S_{n-p}} $.  This functor will have endomorphisms $ Y_{n-p+1}, \dots, Y_n $ and $ (n-p+1 \ n-p+2), \dots, (n-1\ n) $.  It is easily seen that they define an action of a degenerate affine Hecke algebra $\overline{H}_p $ on $ Res^{S_n}_{S_{n-p}} $.  For any $ \mu $ with $ \sum \mu_i = n $, the functor $ E^p_i : Rep(S_n)_\mu \rightarrow Rep(S_{n-p})_{\mu - p \alpha_i} $ is a direct summand of $ Res^{S_n}_{S_{n-p}}$ and thus $ E^p_i $ carries an action of $ \overline{H}_p $.  Theorem 3.16 from \cite{R1} explains how we can convert this to an action of the nil affine Hecke algebra $H_p$ (a similar result was obtained by Brundan-Kleshchev \cite{BK}).  Using this result, we can construct the categorical $ \widehat{\sl}_p$-action.  For more details, see section 5.3.7 of \cite{R1}.

\subsection{Cyclotomic quotients}
There is a natural way to construct categorifications of irreducible representations of $ \g $ using cyclotomic quotients of KLR algebras.  For each dominant weight $ \lambda = \sum n_i \omega_i \in X_+ $ and each $ \nu \in \N I $, let $ R_\nu(\lambda) $ be the quotient of $ R_\nu $ by the ideal generated by all diagrams of the form
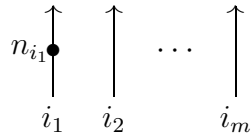
\begin{equation*}
\cyclo
\end{equation*}

Let $ \cV(\lambda)_\mu = R_{\lambda - \mu}(\lambda)$-$\pmod$ be the category of projective modules over the cyclotomic quotients.  Note that there is an action of $\cU^- \g $ on $ \cV(\lambda) $ coming from maps $$ R_{\lambda - \mu}(\lambda) \otimes R_\nu \rightarrow R_{\lambda - (\mu - \nu)}(\lambda) $$ which are given by horizontal concatenation (here $ \cU^- \g $ is defined in the same fashion as $ \cU^+ \g $).  In particular, we have functors $ F_i : \cV(\lambda) \rightarrow \cV(\lambda) $.

The following result was conjectured by Khovanov-Lauda \cite{KL1} and was proved by Kang-Kashiwara \cite{KK,Ka} and Webster \cite{W}.  The conjecture was motivated by the work of Ariki \cite{Ar} who constructed naive categorical $ \widehat{\sl_n} $ actions using cyclotomic Hecke algebras. 
\begin{theo}
The functors $ F_i $ admit biadjoints $ E_i $ and this defines a categorical $ \g$-action on $ \cV(\lambda) $.  Moreover, $ \cV(\lambda) $ categorifies the irreducible representation $ V(\lambda) $ of highest weight $\lambda $.
\end{theo}

Rouquier has proved \cite{R2} that a slight generalization of $ \cV(\lambda) $ is the universal categorical $ \g$-representation with highest weight $\lambda $. Also, Lauda-Vazirani \cite{LV} constructed the crystal of $ V(\lambda) $ using the simple modules over the algebras $ R_\nu(\lambda) $.

\begin{rema}
Webster \cite{W} has generalized this construction.  For any sequence $ \lambda_1, \dots, \lambda_n $, he has constructed certain diagrammatic algebras $ R_\nu(\lambda_1, \dots, \lambda_n) $ whose categories of projective modules admit a categorical $ \g $-action as above.  This construction categorifies the tensor product representation $ V(\lambda_1) \otimes \cdots \otimes V(\lambda_n) $.
\end{rema}

\subsection{Geometric examples}
It is natural to generalize the construction of the categorical $ \sl_2$-action on sheaves on Grassmannians (Theorem \ref{th:sl2actionGrass}).  The generalization uses Nakajima quiver varieties.

For each dominant weight $\lambda $, there exists a (disconnected) Nakajima quiver variety $ Y(\lambda) = \cup_{\mu \in X} Y(\lambda, \mu) $.  It is a moduli space of framed representations of a doubled quiver with preprojective relation.  Nakajima \cite{N} has constructed an action of $ \g $ on $ H^*(Y(\lambda)) $.  This motivated the question of constructing categorical actions of $ \g $ on categories defined using $ Y(\lambda) $.

The variety $ Y(\lambda,\mu) $ is almost a cotangent bundle --- in fact, it can be viewed as an open subset of a cotangent bundle to a certain stack $ M(\lambda, \mu) $.  This motivated Zheng \cite{Z} to define a certain category of constructible sheaves $ D(\lambda, \mu) $ on $ M(\lambda, \mu) $ which should carry a categorical $ \g $ action.

\begin{exem}
In the case $ \g = \sl_2$ and $ \lambda = n $, then $ Y(n, n-2k) = T^* G(k, \C^n) $ and $ D(n, n-2k) = D_c^b(G(k,\C^n))$.
\end{exem}

\begin{theo}
The categories $ D(\lambda, \mu) $ carry a categorical $ \g$-action.
\end{theo}
Zheng \cite{Z} proved a weaker version of this theorem (he only established a naive categorical action).  The above statement was proven by Rouquier \cite{R2} using Theorem \ref{th:VV}.  Webster \cite{W2} also explained that the category $ D(\lambda, \mu) $ can be viewed (under Riemann-Hilbert correspondence and Hamiltonian reduction) as a category of modules over a deformation quantization of $ Y(\lambda, \mu) $.

\end{document}